\documentclass[reqno,10pt]{article}

\usepackage{hyperref}
\usepackage{amsmath}
\usepackage{amssymb}
\usepackage{mathrsfs}

\newtheorem{thm}{Theorem}[section]
\newtheorem{lemma}[thm]{Lemma}

\newtheorem{corol}[thm]{Corollary}
\newtheorem{propos}[thm]{Proposition}
\newtheorem{rema}{Remark}[section]
\def\bp{\begin{propos}}
\def\ep{\end{propos}}
\def\bt{\begin{thm}}
\def\et{\end{thm}}
\def\bco{\begin{corol}}
\def\eco{\end{corol}}
\def\bl{\begin{lemma}}
\def\el{\end{lemma}}
\def\br{\begin{rema}}
\def\er{\end{rema}}
\def\be{\begin{equation}}
\def\ee{\end{equation}}
\def\ba{\begin{array}}
\def\ea{\end{array}}
\def\bena{\begin{eqnarray}}
\def\eena{\end{eqnarray}}

\def\P{{\mathbb P}}
\def\E{{\mathbb E}}
\def\R{{\mathbb R}}
\def\Z{{\mathbb Z}}

\def\1{I}
\def\imath{\textbf{i}}
\def\jmath{\textbf{j}}
\def\chi{\zeta}

\def\a{{\alpha}}

\def\QED{\hfill$\square$\vskip 3mm}

\def\Dp{\displaystyle}

\def\({\left(}
\def\){\right)}

\begin{document}

\title{On The Waiting Time for A M/M/1 Queue\\ with Impatience
\\[5mm]
\footnotetext{*Correspondence author}
\footnotetext{AMS classification (2000): Primary 60K 25, 60F 05; secondary 93A 30} \footnotetext{Key words
and phrases: M/M/1 queue, last in first out, waiting time, tail probability, Laplace transforms
}
\footnotetext{Research supported in part by the Natural Science Foundation of China (under
grants 11271356, 11471222, 11671275) and the Foundation of Beijing Education Bureau (under grant KM201510028002) }}

\author{Feng Wang,\ \  Xian-Yuan Wu$^*$}\vskip 10mm
\date{}
\maketitle

\begin{center}
\begin{minipage}{10.5cm}

\noindent School of Mathematical Sciences, Capital
Normal University, Beijing, 100048, China. Email:
\texttt{fwang@cnu.edu.cn}, \texttt{wuxy@cnu.edu.cn}
\end{minipage}
\end{center}
\vskip 5mm
{\begin{center} \begin{minipage}{10.8cm}
{\bf Abstract}: This paper focuses on the problem of modeling the correspondence pattern for ordinary people. Suppose that letters arrive at a rate $\lambda$ and are answered at a rate $\mu$. Furthermore, we assume that, for a constant $T$, a letter is disregarded when its waiting time exceeds $T$, and the remains are answered in {\it last in first out} order. Let $W_n$ be the waiting time of the $n$-th {\it answered} letter. It is proved that $W_n$ converges weekly to $W_T$, a non-negative random variable which possesses a density with {\it power-law} tail when $\lambda=\mu$ and with exponential tail otherwise. Note that this may provide a reasonable explanation to the phenomenons reported by Oliveira and Barab\'asi in \cite{OB}.


\end{minipage}
\end{center}}

 \vskip 5mm
\section{Introduction and statement of the results}
\renewcommand{\theequation}{1.\arabic{equation}}
\setcounter{equation}{0}

In 2005, Oliveira and Barab\'asi \cite{OB} reported their research results on the correspondence patterns of Darwin and Einstein: during their lifetimes,
Darwin and Einstein answered a fraction of letters they received (the over-all response rate being 0.32 and 0.24, respectively), and the distributions of response times to letters are both well approximated with a power-law tail that has an exponent $\a=3/2$. The classical M/M/1 process \cite[Section 2.2]{GSTH}, which assumes that letters arrive at a rate $\lambda$ and are answered at a rate $\mu$, can be used to model their correspondence patterns. Under some service discipline, the waiting-time density of the M/M/1 process may follow $f(t)\sim t^{-3/2}\exp(-t/t_0)$ for $\lambda\leq \mu$ (see \cite{AW}), which predicts a power-law waiting time for the critical regime $\lambda=\mu$, when $t_0=\infty$. However, by the response rates 0.32 and 0.24, we have $\lambda>\mu$ and this places the model in the supercritical regime, where a finite fraction of letters are never answered. Oliveira and Barab\'asi pointed out in \cite{OB} that numerical simulations indicate that in this supercritical regime the waiting-time distribution of the responded letters also follows a power law with exponent $\a=3/2$.

Clearly, \cite{OB} proposed such questions on modeling the correspondence patterns of human being: how does the ordinary people prioritize the correspondence in need of a response? and does a usual priority principle really lead to a power-law waiting time in the supercritical regime? To partially answer the above questions, or at least, to provide some useful evidence for understanding the above questions, in the present paper, we introduce a special queueing system with a service discipline, which corresponds an usual priority principle of ordinary people, and then study the waiting time for served customers, especially in the supercritical regime.

Now, let's consider the usual M/M/1 queue, the simplest queueing model used in practice. Suppose that the arrivals occur in a Poisson process with rate $\lambda$ and the service times of the unique server have an exponential distribution with parameter $\mu$. We put the service discipline as follows: for some fixed $T>0$, a customer leaves the queue when his waiting time exceeds $T$; the remains are served on the {\it last in first out} principle, namely, the principle serves customers one at a time, and the customer with the shortest waiting time will be served first. Note that our model is a queueing system with impatient customers, which was first proposed by Barrer \cite{B1,B2}, then followed by \cite{F,Ju1,Ju2,KT} {\it etc.}. For details on queues with impatience, one may refer to \cite{WLJ} and the references therein.

Assume that at time $0$, a customer (the $0$-th customer) arrives and the queueing system begins to work.

Let $\{X_n:n\geq 1\}$ be an i.i.d sequence of exponential distribution random variables with parameter $\lambda$, independent of $\{X_n\}$, $\{Y_n:n\geq 1\}$ be another i.i.d sequence of exponential distribution random variables with parameter $\mu$. Let $S^X_n=\Dp\sum_{i=1}^nX_i$ and $S^Y_n=\Dp\sum_{i=1}^n Y_i$. In our system, the $n$-th customer arrives at time $S^X_n$.

Now, fix some $\lambda$, $\mu>0$ and $0<T\leq \infty$. Let $N_t$ be the {\it queue length}, namely, the total number of customers in the system at time $t$. It is clear that $\{N_t:t\geq 0\}$ is no longer a Markov process for any $0<T<\infty$; and in the case of $T=\infty$, $\{N_t\}$ degenerates to the classical birth-death process with birth rate $\lambda$ and death rate $\mu$.
For any $n\geq 0$, Denote by $D_n$ the waiting time of the $n$-th customer, namely, $D_n$ is the usual waiting time when he is finally served before his waiting time exceeds $T$, otherwise, $D_n=\infty$. By the setting given in the last paragraph, one has $N_0=1$ and $D_0=0$. Let
\be\label{1}\tau=\tau(\lambda,\mu,T)=\inf\{n\geq 1: D_n=0\}.\ee
Denote $q_k=\P(\tau=k)$, $k\geq 1$, and let $M=M(\lambda,\mu,T)=\E(\tau)$ be the expectation of $\tau$.

First of all, we have the following proposition on $\E(\tau)$.

\bp\label{p1} i) for any $\lambda,\ \mu>0$ and $0<T<\infty$, we have
\be\label{2}M=M(\lambda,\mu,T)=\E(\tau)<\infty.\ee
ii) for any given $\mu>0$, $M=M(\lambda,\mu,T)$ increases in $\lambda$ and $T$. Furthermore,
for any given $\lambda$ and $\mu$,
\be\label{3}\lim_{T\rightarrow\infty}M(\lambda,\mu,T)=M(\lambda,\mu,\infty)\left\{\ba{ll}<\infty,&\ \rm{if}\ \lambda<\mu;\\[3mm]
=\infty,&\ \rm{if}\ \lambda\geq\mu.\ea\right.\ee\ep


Let $\bar D_n$ be the number of customers who have already been in the queue when the $n$-th customer arrives. Then, almost surely,
$D_n=0$ if and only if $\bar D_n=0$, and $\tau$ is the first return time of state $0$ for process $\{\bar D_n:n\geq 0\}$.

By i) of Proposition~\ref{p1} and a renewal limit theorem for general stochastic processes (see \cite{Be} or \cite[Section 8, Chapter 11]{Fe}), the limit \be\label{4'}\lim_{n\rightarrow\infty}\P(D_n=0)=\lim_{n\rightarrow\infty}\P(\bar D_n=0) \ \rm{exists}\ee and will be proved to be $1/M$ ( see Lemma~\ref{l1} in Section 2).

We point out that, for any $0<T<\infty$, by i) of Proposition~\ref{p1} and the renewal limit theorem in \cite{Be}, our system tends to {\it statistical equilibrium}, namely the queue length $N_t$ converges in distribution to a $\Z_+$-valued random variable $N$ as $t\rightarrow\infty$.

Let $I_1$ be the modified Bessel function of first kind for real and positive $t$, see \cite{O}, given by
$$I_1(t)=\sum_{m=0}^{\infty}\frac 1{m!(m+1)!}\(\frac t 2\)^{2m+1}\sim \frac{e^{t}}{\sqrt {2\pi t}}, \ \rm{as}\ t\rightarrow\infty.$$ Note that for two functions $a(t)$ and $b(t)$, $a(t)\sim b(t),\ t\rightarrow\infty $ means that $\lim_{t\rightarrow\infty} a(t)/b(t)=1$.

Now, we can state our main theorem as follows.

\bt\label{th1} For any $\lambda,\ \mu>0$ and $0<T\leq\infty$, $W_n$, the waiting time of the $n$-th served customer, converges in distribution to a non-negative random variable $W_T$. Furthermore, the distribution function $F_T$ of $W_T$ satisfies $F_T(x)=0$, for $x<0$, $F_T(x)=1$ for $x> T$, and
\be\label{5}F_T(x)=\frac 1{C(T)}\left[\frac 1M+\frac 1{\rho\vee 1}\(1-\frac 1M\)\int^x_0 f_\rho (t)dt\right],\ee for $0\leq x\leq T$, where $C(T)$ is the normalization constant, $\rho=\lambda/\mu$ and
\be\label{6}f_\rho(t)=\sqrt{\rho\vee\rho^{-1}}\ \frac 1t\ e^{-(\lambda+\mu)t}I_1(2t\sqrt{\lambda\mu}),\ \ t>0.\ee\et

\br The function $f_\rho$ given in (\ref{6}) is the probability density of $D$, the length of the busy period in the classical M/M/1 system when $\rho\leq 1$, and is the {\bf conditional} probability density of $D$ conditioned on $D<\infty$ when $\rho>1$. Note that in case of $\rho>1$, $\P(D<\infty)=\rho^{-1}$.  \er

\br For any large $0<T\leq \infty$, $f_T(t)$, the density of $W_T$ has a power-law tail with exponent $\alpha=3/2$ in the critical regime $\lambda=\mu$ as $t\uparrow T$. In both subcritical and supercritical cases, $f_T(t)$ decays exponentially fast.   \er

But, what are the correspondence patterns of Darwin and Einstein? Numerical simulations indicate that they always keep themselves in the critical regime. A reasonable explanation may be the following: letters arrive according to a Poisson process with rate $\lambda=\lambda_1+\lambda_2>\mu$, where $\lambda_1<\mu$ is the rate of letters heard from friends and family members, $\lambda_2$ is the rate of letters heard from the strangers. As the most distinguished scientists in their research fields, Darwin and Einstein received too many letters from the strangers, so they had to {\bf ignore} such a received letter with probability $1-(\mu-\lambda_1)/\lambda_2$, such that, in their eyes, letters arrived according to a Poisson process with rate $\lambda_1+\lambda_2\times {(\mu-\lambda_1)}/{\lambda_2}=\mu$.

\vskip 5mm
\section{Proofs}
\renewcommand{\theequation}{2.\arabic{equation}}
\setcounter{equation}{0}

In this section we will prove Proposition~\ref{p1} and Theorem~\ref{th1}. Before we give a proof to Proposition~\ref{p1}, we introduce an algorithm to obtain $\tau$.

Given $\lambda, \mu$ and $T$, let $S^X_n,\ S^Y_n$ be defined as in Section 1 and let $\bar S^X_n=S^X_n+T$. For any integer $k\geq 1$, define the random index set $N_k$ as
$$N_k:=\left\{n: S^Y_{k-1}\leq S^X_n<S^Y_k\right\},\ \ \rm{where\ }\ S^Y_0=0.$$
Let $\bar N_k=\cup_{i=1}^k N_i$ and write $n_k$ as the largest element in $\bar N_k$ when $\bar N_k\not=\emptyset$.

Now we will start a procedure to define a sequence of index sets $\{J_k:k\geq 1\}$, we write $j_k$ as the largest element of $J_k$ when $J_k\not=\emptyset$.
Let $J_1=\{n\in N_1:\bar S^X_n>S^Y_1\}$. If $J_1\not=\emptyset$, then let $J'_1=J_1\setminus\{j_1\}$ and define $J_2=\{n\in J'_1\cup N_2: \bar S^X_n>S^Y_2\}$; otherwise, the procedure is stopped. By induction, for any $k\geq 2$, if $J_k$ has already been defined and $J_k\not=\emptyset$, then let $J'_k=J_k\setminus\{j_k\}$ and define $J_{k+1}=\{n\in J'_k\cup N_{k+1}:\bar S^X_n>S^Y_{k+1}\}$; otherwise, the procedure is stopped. Let
\be\label{7}\tau_1=\inf\{k\geq 1:J_k=\emptyset\},\ \ \rm{(set\ }\inf\emptyset=\infty\rm{ )},\ee then $$\tau=\left\{\ba{lll}&\hskip-4mm n_{\tau_1}+1,&\rm{if}\ \ \tau_1<\infty;\\[2mm]
&\hskip-4mm\infty,&\rm{otherwise}.\ea\right.$$

{\it Proof of i) of Proposition~\ref{p1}.} Suppose that $0<T<\infty$. For any $k\geq 1$, let $A_k$ be the event that $Y_k\geq T$ and $N_k=\emptyset$. Recalling that the random point set $\{S^X_n:n\geq 1\}$ forms a one dimensional Poisson point process in $\R_+$, we know that $\{A_k:k\geq 1\}$ are mutually independent and, for any $k\geq 1$,
$$p_0:=\P(A_k)=\Dp\int_T^{\infty}\P(N_k=\emptyset\mid Y_k=t)\mu e^{-\mu t}dt=\int_T^{\infty}e^{-\lambda t}\mu e^{-\mu t}dt=\Dp\frac{\mu}{\lambda+\mu} e^{-(\lambda+\mu) T}.$$
Let $$\tau_2=\inf\{k\geq 1:A_k \ \rm{occurs}\}.$$ First, we point out that $\tau_2$ is a geometry distribution random variable with parameter $p_0$,
hence, we have
\be\label{8}\E(\tau_2)=\frac 1{p_0}=\frac\mu{\lambda+\mu}e^{(\lambda+\mu)T}.\ee
Second, by (\ref{7}), we have $\tau_1\leq \tau_2$ and $n_{\tau_1}\leq n_{\tau_2}$. Denote by $h_k$ the cardinality of $N_k$, then $$n_{\tau_2}=\sum_{k=1}^{\tau_2}h_k,$$and \be\label{9}\ba{ll}&\E(n_{\tau_2})=\Dp\sum_{n=1}^{\infty}\E(n_{\tau_2}\mid\tau_2=n)\P(\tau_2=n)\\&=\Dp\sum_{n=1}^\infty\sum_{k=1}^n\E(h_k\mid \tau_2=n)\P(\tau_2=n). \ea\ee
Noticing that $\{\tau_2=n\}=\(\cap_{k=1}^{n-1}A_k^c\)\cap A_n$, we have $\E(h_n\mid \tau_2=n)=0$ and, for any $1\leq k\leq n-1$,
\be\label{10}\ba{ll}\E(h_k\mid \tau_2=n)&=\E(h_k\mid A_k^c)=\E(h_k\mid h_k>0)\P(h_k>0\mid A_k^c)\\[3mm]&\leq \E(h_k\mid h_k>0)=\Dp\frac{\E(h_k)}{\P(h_k>0)},\ea\ee where the first equality comes from the fact that $h_k$ is independent of $A_{k'}$ for any $k'\not=k$. Clearly,
$$\E(h_k)=\Dp\int_{0}^\infty\E(h_k\mid Y_k=t)\mu e^{-\mu t}dt=\int_{0}^\infty \lambda t\mu e^{-\mu t}dt=\frac \lambda\mu; \ \rm{and}$$
$$\P(h_k>0)=\Dp\int_0^\infty\P(h_k>0\mid Y_k=t)\mu e^{-\mu t}dt=\int_0^\infty\(1-e^{-\lambda t}\)\mu e^{-\mu t}dt=\frac\lambda{\lambda+\mu}.$$
Thus, by (\ref{9}), (\ref{10}) and then by (\ref{8}),
$$\ba{ll}\E(\tau)&=\Dp\E(n_{\tau_1})+1\leq\E(n_{\tau_2})+1\leq \frac{\lambda+\mu}\mu\sum_{n=1}^\infty(n-1)\P(\tau_2=n)+1\\[5mm]&\leq \Dp \frac{\lambda+\mu}\mu\E(\tau_2)=e^{(\lambda+\mu)T}<\infty. \ea$$ \QED

{\it Proof of ii) of Proposition~\ref{p1}.} Given $\mu$ and $T$. To see $M(\lambda,\mu,T)$ increases in $\lambda$, one only need to notice the fact that two Poisson processes with parameters $\lambda_1$ and $\lambda_2$, $\lambda_1<\lambda_2$, respectively can be coupled together in such a way that all arrivals of the former forms a part arrivals of the later. Note that this follows from the {\it infinitely divisible} property of Poisson distribution \cite{JKK}.  For details on coupling technique in probability theory, please refer to \cite{L}. Then, the corresponding monotonicity follows from (\ref{7}), the definition of $\tau_1$.

For given $\lambda$ and $\mu$, the fact that $M(\lambda,\mu,T)$ increases in $T$ follows directly from the definition of $\tau_1$.

 Now we are in the way to prove (\ref{3}), note that the $\lambda\geq \mu $ part of (\ref{3}) only need to be proved for $\lambda=\mu$, the case of $\lambda>\mu$ follows from the monotonicity in $\lambda$ proved above.

We first declare that $M(\lambda,\mu,\infty)<\infty$ for $\lambda<\mu$ and $=\infty$ for $\lambda=\mu$. In fact, in case of $T=\infty$, our system degenerates to the classical M/M/1 queue. The standard argument of the birth-death process (see \cite{A,CH}) tells us that \be\label{11}\P(\tau(\lambda,\mu,\infty)<\infty)=1,\ \rm{for}\  \lambda\leq\mu,\ee and $$M(\lambda,\mu,\infty)=\E(\tau(\lambda,\mu,\infty))\left\{\ba{ll}<\infty,& \rm{for}\ \lambda<\mu;\\[2mm] \infty, & \rm{for}\ \lambda=\mu.\ea\right.$$

Second, by (\ref{11}) and the algorithm we have used to obtain $\tau$, we have $\tau(\lambda,\mu,T)$ increases in $T$, and
$$\lim_{T\rightarrow\infty}\tau(\lambda,\mu,T)=\tau(\lambda,\mu,\infty), \ \ \rm{a.s.}.$$ Thus, (\ref{3}) follows from the monotone convergence theorem.\QED

Before we give a proof to Theorem~\ref{th1}, we have to prove the following lemma. Note that this lemma plays the key role in the proof of  Theorem~\ref{th1}.

\bl\label{l1}Suppose that $\lambda,\mu>0$ and $0<T\leq\infty$, then
\be\label{12}\lim_{n\rightarrow\infty}\P(D_n=0)=\frac 1{M(\lambda,\mu,T)}.\ee\el

{\it Proof.} In case of $T=\infty$, we are dealing with a classical M/M/1 system, (\ref{12}) follows from a standard argument for birth-death process \cite{A,CH}. Now, we suppose
$0<T<\infty$.

Let $P^{(n)}_{00}=\P(D_n=0),\ n\geq 0$. Then $P^{(0)}_{00}=1$, and for any $n\geq 1$,
\be\label{13}P^{(n)}_{00}=\sum_{k=1}^nq_kP^{(n-k)}_{00},\ee where $q_k=\P(\tau=k)$. (\ref{13}) indicates that the sequence $\{P^{(n)}_{00}:n\geq 0\}$ is iteratively determined by $\{q_k:k\geq 1\}$ and its initial value $P^{(0)}_{00}=1$.

Provided (\ref{4'}) and (\ref{13}), the lemma may follows from a standard argument on generation function and the Abel's Theorem (see \cite[page 12]{CM}). Here, we prefer to give it a probabilistic proof, one will see that our proof mainly depends on (\ref{13}), and (\ref{4'}) is only a consequence of it.

By the basic theory on discrete-time Markov chains \cite[Section 1.2]{CM}, to prove (\ref{12}), it suffices to construct a $\Z_+$-valued discrete-time Markov chain $\{\xi_n:n\geq 0\}$ such that $\{\xi_n\}$ is {\it ergodic} and
\be\label{14}f_{00}^{(k)}:=\P(\tau^+_0=k\mid \xi_0=0)=q_k,\ \ k\geq 1,\ee
where $\tau^+_0=\inf\{n\geq 1:\xi_n=0\}$ is the first return time of state $0$. Note that {\it ergodic} means {\it irreducible, aperiodic } and {\it positive recurrent} as usual.

To this end, we give the following transition matrix $P=(P_{ij})$ to $\{\xi_n\}$:
$$P_{ij}=\left\{\ba{ll}\Dp\frac{q_{i+1}}{1-\sum_{k=1}^iq_k},&j=0;\\[3mm]
1-P_{i0},&j=i+1;\\
0,&\rm{else}.\ea\right.$$
It is straightforward to check that $\{\xi_n\}$ with the above transition matrix $P$ is ergodic and satisfies (\ref{14}).
Thus we finish the proof of the lemma.\QED

{\it Proof of Theorem~\ref{th1}.} For any $0\leq x\leq T$, one has
\be\label{17} \P(D_n\leq x)=\P(D_n=0)+\P(D_n\leq x\mid D_n>0)\P(D_n>0).\ee
First of all, $D_n>0$ if and only if, at $S^X_n$, the arrival time of the $n$-th customer, the unique server is occupied, namely, $\bar D_n>0$. It is clear that, on the {\it last in first out} principle, if $\bar D_n>0$, then $D_n$ does not depend on the exact value of $\bar D_n$. Second, for any $0<x\leq T$, $D_n\leq x$ means that all customers arrived after, but were served before the $n$-th one are finally served in time $x$ ($\leq T$). In other word, the $n$-th customer finally gets into the sever after a whole {\it busy period} finishes in time $x$.
Hence, by the {\it memoryless property} of the exponential distribution, for any $0<x\leq T$ and for any $n\geq 1$,
$$\P(D_n\leq x\mid D_n>0)=\P(D\leq x)=:F_{D}(x),$$
where $D$ is the length of the busy period of the classical M/M/1 system. Note that, in our setting, $D=S^Y_{\tau_3}$, with
\be\label{4}\tau_3:=\inf\left\{n:S^Y_n<S^X_n\right\}.\ee As calculated in \cite{Xu},
\be\label{20}F_D(x)=\sum_{m=1}^\infty\int_0^x e^{-\lambda t}\frac{(\lambda t)^{m-1}}{m!}dF_{S^Y_m}(t)=\sum_{m=0}^\infty\int_0^x \frac{\mu(t\sqrt{\lambda\mu} )^{2m}}{m!(m+1)!}e^{-(\lambda+\mu) t}dt.\ee Where $F_{S^Y_m}$ is the distribution function of $S^Y_m$.

In case of $\lambda\leq \mu$, by (\ref{4}), we have $\P(D<\infty)=\P(\tau_3<\infty)=1$, so $F_D$ is a probability distribution function and
\be\label{15'} f_{D}(t):=F_{D}'(t)=\sqrt{\frac \mu\lambda}\ \frac 1t\ e^{-(\lambda+\mu)t}I_1(2t\sqrt{\lambda\mu})\ee is a probability density.

To finish the proof of the theorem, we will introduce $\Gamma$, the Laplace transform of $F_{D}$ defined by the following Lebesgue-Stieltjes integration
$$\Gamma(s)=\int^\infty_0 e^{-st}dF_{D}(t), \ \ \rm{Re}(s)>0,$$ where $\rm{Re}(s)$ is the real part of the complex number $s$.

It is calculated directly that (see \cite[Eq. (51)]{Xu}),
\be\label{15}\Gamma(s)=\frac 1{2\lambda}\left[\lambda+\mu+s-\sqrt{(\lambda+\mu+s)^2-4\lambda\mu}\right].\ee Clearly, $\Gamma$ can be uniquely inverted to give the probability distribution function $F_D$ in case of $\lambda\leq\mu$.

In case of $\lambda>\mu$, by symmetry and the linearity of the (inverse) Laplace Transform, the following
$$\tilde\Gamma(s)=\rho\Gamma(s)=\frac 1{2\mu}\left[\lambda+\mu+s-\sqrt{(\lambda+\mu+s)^2-4\lambda\mu}\right]$$
can be inverted to give the following probability distribution function
\be\label{16}\tilde F_D(x):=\rho F_D(x)=\P(D\leq x)/\P(D<\infty)=\P(D\leq x\mid D<\infty),\ee
and $\tilde f_D=\rho f_D$ is the probability density corresponding to $\tilde F_D$.

By (\ref{16}), one has $\P(D<\infty)=\rho^{-1}$ and $\tilde F_D$ is the conditional distribution function of $D$ conditioned on $D<\infty$.

Write $f_D$ and $\tilde f_D$ in the unified form $f_\rho $, which is given in (\ref{6}). Then, by (\ref{17}), the conditional distribution of $D_n$ conditioned on $D_n\leq T$ is
\be\label{18}\ba{rl} F_n(x)&\hskip-4mm:=\P(D_n\leq x\mid D_n\leq T)\\[3mm]&\hskip-4mm=\frac 1{C_n(T)}\left[\P(D_n=0)+\P(D_n>0)\frac 1{\rho\vee 1}\int^x_0 f_\rho (t)dt\right],\ea\ee for $0\leq x\leq T$, where
$C_n(T)$ is the normalization constant. The theorem follows immediately form (\ref{18}) and Lemma~\ref{l1}.
\QED


\vskip10mm

\end{document}